\documentclass[12]{article}
\usepackage{amsmath}
\usepackage{amsfonts}
\usepackage{amssymb}
\usepackage{amsthm}
\usepackage{color}

\newtheorem{ex}{Example}

\newtheorem{df}{Definition}[section]

\newcommand{\R}{{\rm I}\kern-0.18em{\rm R}}
\newcommand{\1}{{\rm 1}\kern-0.25em{\rm I}}
\newcommand{\E}{{\rm I}\kern-0.18em{\rm E}}
\newcommand{\p}{{\rm I}\kern-0.18em{\rm P}}


\usepackage{epsfig}

\usepackage{fancyhdr}

\usepackage{graphics}

\usepackage{amsopn}  

\usepackage{amssymb}
\usepackage{amsxtra}
\usepackage{amsfonts}

\usepackage{afterpage}
\usepackage{eucal}
\usepackage{eufrak}

\usepackage{enumerate}

\allowbreak

\def\fnote#1{\footnote}

\newcommand{\bea}{\begin{eqnarray}}

\newcommand{\eea}{\end{eqnarray}}

\newcommand{\beas}{\begin{eqnarray*}}
\newcommand{\eeas}{\end{eqnarray*}}

\author{Lev B Klebanov\footnote{Department of Probability and Statistics, MFF, Charles University, Prague-8, 18675, Czech Republic},  Lenka Sl\'amov\'a\footnote{Department of Probability and Statistics, MFF, Charles University, Prague-8, 18675, Czech Republic and Stony Brook University, New York, U.S.A.}, Ashot Kakosyan \footnote{Yerevan State University, Yerevan, Armenia} , Gregory Temnov \footnote{School of Mathematical Science, University College Cork, Cork, Ireland}}
\title{Casual Stability of Some Systems of Random Variables}
\date{}

\begin{document}
\maketitle

\begin{abstract}

Self-similarity of systems is very popular and intensively developing field during last decades. To this field belong so-called stable distributions and their generalization. In  \cite{KS} was given an approach to define additive systems with the property of random self-similarity - casual stability (c.s.). Here we continue study the notion of casual stability for additive  systems of random variables (r.v.). We also give a modification of this definition and  spread them on multiplicative  systems of r.v.  and on the system with operations of taking minimum or maximum of r.v. The case of systems with a random number of elements is also considered.
\end{abstract}

{\bf keywords:} stable distributions; casual stable distributions; casual self-similarity

\section{Definitions of casual stability for additive systems}

Everywhere in this paper we will consider r.v. given either on whole real line $\R$, or on non-negative
semi-axes $\R_{+}$. Although it is possible to give some definitions of c.s. for r.v. on $d-$dimensional Euclidean space $\R^d$ or just on Banach spaces, we do not consider such cases here. We use the term casual stability as generalization of the notion of strictly stable r.v. (not stable r.v.). So, sometimes we use term ``stability" instead of ``strictly stability" and hope it will be clear what are we talking about.

As it was mentioned in Abstract, the definition of c.s. for additive systems of random variables was given in \cite{KS}. Let us give it here for the set of r.v. on $\R$ with additive operation for independent r.v.

Let $X$ be a r.v. with cumulative distribution function (c.d.f.) $F(x)$. Characteristic function (ch.f.) $f(t)$ of $X$ may be written as
\[  f(t)= \int_{-\infty}^{\infty}e^{i t x} dF(x)= \int_{0}^{\infty}(e^{-i t})^{|x|} d(1-F(-x))+ \int_{0}^{\infty}(e^{i t})^{x} dF(x).\]
 Now we make ``random normalization" in the following way. We change $e^{i t}$ in the second integral by a ch. f. $g(t)$, but it the first integral we change $e^{-i t}$ by $g(-t)$. The function obtained denote by ${\tilde f}(t)$ and say it is $g$-normalization of $f(t)$.
 
 \begin{df}\label{Def1} Let $\{g_n (t)\;n =1,2, \ldots\}$ be a sequence of ch.f. We say, a ch.f. $f(t)$ is c.s.  ch.f. with respect to the sequence $\{g_n (t)\;n =1,2, \ldots\}$ (or just c.s.) if for any integer
 $n \geq 1$ 
\begin{equation}\label{eq1}
  {\tilde f}^n_n(t) = f(t), \; t \in \R , 
  \end{equation}
 where $ {\tilde f}_n(t)$ is $g_n$-normalization of $f(t)$.
\end{df}

From Definition \ref{Def1} it follows that {\it casual stable distributions are infinitely divisible}.

Let us note, that the normalization used in definition of discrete stable distribution, is a particular case of general random normalization. As far as we know, the definition of discrete stability and  corresponding normalization were given in \cite{StvHarn79}; see also \cite{ChSch98a}.

The Definition \ref{Def1} may be reformulated in terms of r.v. Namely, let $X$ be a r.v. having ch.f. $f(t)$, and ${\tilde X}$ be r.v. with ch.f. ${\tilde f}(t)$- $g$-normalization of $f(t)$. We say $X$ is c.s. r.v. if for any integer $n \geq 1$ the following equality in distribution 
\begin{equation}\label{eq2}
 X \stackrel{d}{=} \sum_{j=1}^n {\tilde X}_j(n) 
 \end{equation}
holds. Here $ {\tilde X}_j(n)$, $j=1, \ldots , n$ are independent identically distributed (i.i.d.) random variables with ch.f. ${\tilde f}_n$.

Of course (and it was mentioned in \cite{KS}), classical strictly stable distributions are c.s. too.

The Definition  \ref{Def1} gives only one possible way to define c.s.r.v. The matter is such that classical strictly stable r.v. may be defined in many equivalent  ways. However, such ways lead to different classes of c.s. r.v. Let us give another definition for a new class of c.s. r.v. We call such r.v. pursuit casual stable (p.c.s.).

\begin{df}\label{Def2}  Let $\{g_n (t)\;n =1,2, \ldots\}$ be a sequence of ch.f. We say, a ch.f. $f(t)$ is p.c.s.  ch.f. with respect to the sequence $\{g_n (t)\;n =1,2, \ldots\}$ (or just p.c.s.) if for any integer
 $n \geq 1$ 
\begin{equation}\label{eq3}
  {\tilde f}_n(t) = f(t)^n, \; t \in \R , 
  \end{equation}
 where $ {\tilde f}_n(t)$ is $g_n$-normalization of $f(t)$.
\end{df}

Definition \ref{Def2} is very similar with Definition \ref{Def1}. The only difference between them is that the power $n$ in Definition \ref{Def1} is applied to $g_n$-normalization of $f$, while in Definition Definition \ref{Def2} the power is applied to $f$ itself. Of course, for classical case the definitions are equivalent, because to obtain the sequence for second definition it is sufficient to take the sequence of inverse values of the sequence from the first definition, and vice versa.  Later we will see, that for general case these definitions are essentially different. Particularly, {\it pursuit casual stable distribution may be not infinite divisible}. Corresponding example will given in corresponding section below.  Of course, it is possible to use as a starting point other definitions of strictly stable r.v., but we will not consider them in this paper.

\section{Casual stability for multiplicative type of systems}
\setcounter{equation}{0}

From the Definition \ref{Def1} (especially, from equations (\ref{eq1}) and (\ref{eq2})) it follows that the definition of c.s. r.v. is based on summation of independent r.v. Moreover, we  used ch.f. (and \lq\lq{}elements\rq\rq{} $e^{\pm i t}$ in them) as a tool for work with sums of i.i.d. r.v.  Of course, some other operations on r.v. are of essential interest. Below we give definitions of c.s. r.v. for such operations as multiplication of independent r.v., and that of taking minimum and maximum of r.v. Obviously, instead of ch.f. it is naturally to use Mellin transformation of c.d.f. for multiplication, and survival function and c.d.f. for the cases of minimum and maximum, correspondingly.

Let us start with definition of multiplicative c.s. r.v. taking values on $\R_{+}$. Of course, the main problem is to define an analogue of $g$-normalization for this case. It is well-known that Mellin transform in probability theory is defined for this case as
\begin{equation}\label{eq4}
M(u)=\E\{X^u\}=\int_0^{\infty}x^u dF(x),
\end{equation}
where $F(x)$ is c.d.f. of r.v. $X$, and $u$ is complex number such that the integral in (\ref{eq3}) converges. This transform seems to be applicable for studying products of independent r.v. because Mellin transform of the product of two independent r.v. equals to product of their Mellin transforms.  It is obvious, that 
\[  M(u)=\E \exp(u \log(X)) , \]
and $M(u)$ coincides for $u=it$ with ch.f. of r.v. $\log(X)$. Therefore we may just apply Definitions \ref{Def1} and \ref{Def2} to r.v. $\log(X)$. However, we give  below the definitions in terms of Mellin transform, because it gives ideas for further generalizations. 

Let $X>0$ be a r.v. with Mellin transform
\[ M(u)=\int_0^{\infty}x^u dF(x)=\int_0^{1}x^u dF(x)+\int_1^{\infty}x^u dF(x)= \]
\[ = \int_0^{1}(e^{-u})^{\ln(1/x)} dF(x)+\int_1^{\infty}(e^u)^{\ln(x)} dF(x) =M_1(u)+M_2(u). \]
Following \cite{KS}, we interpret $e^u$ as Laplace transform of degenerate distribution concentrated at point $1$ (It would be more natural to use $\xi = \ln u$ instead of $u$, and speak about Mellin transforms instead of Laplace, but it will give equivalent description). Also we consider a point as a particle, which transforms itself into a flow of particles with a random number of them with Laplace transform $N(u)$. So, we have transformation from $e^u$ to $N(u)$ for $u>0$. Naturally, we define this transformation for negative $u$ in the way that  $e^{-u}$  transforms to $N(-u)$. So, the normalizing transformations is mapping $M(u)$ into ${\tilde M}(u)=M_1(-\ln(N(-u))+M_2(ln(N(u)))$.
Now it is clear that mentioned transform represent N-normalization for the case of multiplicative operation on r.v. Now we are ready to give definition of multiplicative (product) c.s. r.v.

\begin{df}\label{Def3} Let $\{N_n(u),\; n=1,2, \ldots \}$ be a family of Laplace transforms. Suppose that $X>0$ is a r.v. with Mellin transform $M(u)$. We say that $X$ is product c.s. (p.c.s.) if for any integer $n \geq 1$
the following relation 
\begin{equation}\label{eq5}
{\tilde M}_n(u) = M^{1/n}(u),
\end{equation}
where ${\tilde M}_n(u) = M_1(-\ln(N_n(-u)))+M_2(\ln(N_n(u)))$.
holds.
\end{df}

Very similar to Definition \ref{Def2} is the following definition of product (or  multiplicative) pursuit c.s.

\begin{df}\label{Def4} Let $\{N_n(u),\; n=1,2, \ldots \}$ be a family of Laplace transforms. Suppose that $X>0$ is a r.v. with Mellin transform $M(u)$. We say that $X$ is pursuit product c.s. (p.p.c.s.) if for any integer $n \geq 1$
the following relation 
\begin{equation}\label{eq6}
{\tilde M}_n(u) = M^{n}(u),
\end{equation}
where ${\tilde M}_n(u) = M_1(-\ln(N_n(-u)))+M_2(\ln(N_n(u)))$.
holds.
\end{df}

The definitions given in this section are clearly analogous to corresponding definitions for that of additive case. However, it is not enough to consider r.v. with values in $\R_{+}$; we need to propose similar definitions for the case of $\R$. For this aim let us remember the definition of Mellin transform for the case of r.v. $X$ with values in $\R$. Namely,
\begin{equation}\label{eq7}
M_X(u)=\int_0^{\infty}x^u dF_{X_{+}}(x)+\gamma \int_0^{\infty}x^u dF_{X_{-}}(x)=M_{+}(u)+\gamma M_{-}(u),
\end{equation}
where $F_{X_{+}}$ and $F_{X_{-}}$ are c.d.f. of positive and negative parts $X_{+}$ and $X_{-}$ of $X$, and $\gamma$ is a complex number, $\gamma^2=1$. Now Definitions \ref{Def3} and \ref{Def4} may be generalized in the form
\begin{equation}\label{eq8}
{\tilde M}_{+}(u) + {\tilde M}_{-}(u) =M^{1/n}(u).
\end{equation}
and 
\begin{equation}\label{eq9}
{\tilde M}_{+}(u) + {\tilde M}_{-}(u) =M^{n}(u).
\end{equation}
We keep the name p.c.s. for the first case, and p.p.c.s. for the second.

\section{Casual stability for minimum and maximum types of systems}
\setcounter{equation}{0}

Here we give definitions of casual stability and pursuit casual stability for the systems of positive r.v. with operation of taking minimum (maximum) of r.v. We consider only one of these operations, because another may be obtain by passing, say, to inverse values of r.v. The role of characteristic function plays now survival function
\[   {\bar F}(x)= 1- F(x) = \int_0^{\infty}I_{[1,\infty)}(y/x)dF(y),       \]
where $F(x)$ is c.d.f. of a r.v. $X$, and $I_{[1,\infty)}(y)$ is indicator function of the interval $[1,\infty)$.  This is so, because
\[ {\bar F}_X(x) {\bar F}_Y(x) = {\bar F}_{\max(X,Y)}(x), \]
and, particularly, 
\[ I_{(x,\infty)}(y) I_{(z,\infty)}(y) = I_{(\max(x,z),\infty)}(y) \]
what corresponds to degenerate distributions, concentrated at points $x$ and $z$.

Because $I_{[1,\infty)}(y)$ is a c.d.f. of degenerate distribution, it is natural to change it for another (possible, non-degenerate) c.d.f. It is in line of our approach to random normalization.
So, if $X$ is a r.v. with c.d.f. $F(x)$, and $G(x)$ is other c.d.f., then $G$-normalization of $X$ is  r.v. Y with
survival function
\begin{equation}\label{10}
{\tilde F}(x) = {\bar F}(-\ln {\bar G}(x) ).
\end{equation}

Now we can give definitions of minimum casual stable (min-c.s) distribution and pursuit minimum casual stable (p.min-c.s.) distribution.

\begin{df}\label{Def5} Let $\{G_n(u),\; n=1,2, \ldots \}$ be a family of c.d.f. Suppose that $X>0$ is a r.v. with c.d.f $F(u)$. We say that $X$ is min-c.s. if for any integer $n \geq 1$
the following relation 
\begin{equation}\label{eq11}
{\tilde F}_n(u) = F^{1/n}(u),
\end{equation}
(where ${\tilde F}_n(u) = {\bar F}(-\ln {\bar G_n}(x) ) $)
holds.
\end{df}

Very similar to Definition \ref{Def5} is the following definition of pursuit min-c.s.

\begin{df}\label{Def6} Let $\{G_n(u),\; n=1,2, \ldots \}$ be a family of c.d.f. Suppose that $X>0$ is a r.v. with c.d.f $F(u)$. We say that $X$ is min-pursuit c.s. (min-p.c.s.) if for any integer $n \geq 1$
the following relation 
\begin{equation}\label{eq12}
{\tilde F}_n(u) = F^{n}(u),
\end{equation}
(where  ${\tilde F}_n(u) = {\bar F}(-\ln {\bar G_n}(x) ) $)
holds.
\end{df}

We shall not consider minimums of r.v. taking values in $\R$ in this paper.

\section{Systems with random number of elements}

In previously given definitions it was possible to suppose that the system on $n$ randomly normalized elements was equivalent in distribution to one initial (non-normalized) element, or, in the case of pursuit stability, one normalized element was equivalent to the system of $n$ initial elements. Very natural question is the following. What happens if the number $n$ is not fixed, but random. Such question was studied in the paper \cite{KMM84}, where were given definitions of geometrically stable and geometrically infinite divisible r.v. A lot of publications were connected to the study of geometrically stable distributions. We will not give corresponding citations here; the reader may find them, for example, in the book \cite{KKR6}. More general types of random number of elements were studied in \cite{KMM85}, and, in more general situation, in \cite{KlRa} and \cite{Bun}.  In all publication on this problematics only classical normalization (that is, normalization by constant multipliers) was used. However, in some publications were given discrete analogues of some geometric stable distributions (especially, Linnik distributions). For such generalization was used the same normalization as for discrete stable distributions; see \cite{ChSch98b} and \cite{ChSch1}.

Here we propose to use the same definitions as given in previous sections, but use instead of fixed number of elements $n$ a family random variables $\{\nu_p, \; p \in \Delta\}$, as it was proposed in \cite{KlRa}
(see also \cite{KKR6}).  To stress attention on the fact that the number of elements is random, we propose to write $\nu$-casual stable, or similar for other cases.

\section{Examples of casual stable distributions for different types of systems}
\setcounter{equation}{0}

Here we give examples of c.s. distributions for systems, described above. Some of the examples are illustrative only, while other are of essential interest. 

\subsection{Additive systems}

\begin{ex}\label{Ex1} The first example is given by choice of  the sequence of Laplace transforms $g_n(s)=\exp(-a_n s)$. In this case we have just an ordinary normalization, and corresponding casual stable distributions coincide with ordinary positive stable distributions skewed to the right with the parameter $\alpha \in (0,1)$.
\end{ex}

\begin{ex}\label{Ex2} The second example is more interesting from our point of view.  
Let us consider positive stable r.v. with the index of stability $\alpha=1/m$, $m \in N$, move its Laplace transform on $h$ units to the right and make normalization on corresponding measure. We will have Laplace transform of corresponding tempered stable r.v. as:
\[ L(s) =\exp\{ -\lambda^{\alpha}(1+\tan \frac{\pi \alpha}{2})((s+h)^{\alpha}-h^{\alpha})\}. \]
The distribution with this Laplace transform appears to be causal stable with the function
\[ g_n(s) = \exp(h-(1/n (s+h)^{\alpha}+ (n-1)/n h^{\alpha})^{1/\alpha})\]
(it is possible to verify that $g_n$ is Laplace transform of a d.f. in the case when $1/\alpha \in N$). For $h \to 0$ we obtain a classical case of normalization: degenerated distribution at point $1/n^{1/\alpha}$. As a particular case we find that Inverse Gaussian distribution is casual stable too ($\alpha = 1/2$).
\end{ex}

Let us note that the Example \ref{Ex2} was given in \cite{KS}. For $\alpha = 1/2$ we may give an interpretation. For $h=0$ we have L\'evy distribution, which is stable in classical sense. This distribution may be interpreted as a distribution of the moment of the first passage time of Brownian motion for a given level if drift is zero. Similarly, Inverse Gaussian distribution may be interpreted as a distribution of the moment of the first passage time of Brownian motion for a given level if drift is not zero. But this distribution is not stable in classical sense. However, it is casual stable, with distribution of normalization depending on drift. As it was mentioned above, this normalization tends to classical as drift tends to zero.

\begin{ex}\label{Ex3} Let us consider Laplace distribution with ch.f.
\[  f(t)= \frac{1}{1+a^2 t^2} .\]
It may be written as
\[  f(t)= \frac{1}{2}\frac{1}{1-i a t}+\frac{1}{2}\frac{1}{1+i a t}. \]
Make now change $it$ by  $\ln g_n(t)$ in the first summand, and $-it$ by $\ln g_n(-t)$ in the second, where $g_n(t)$ is a ch.f. Supposing that $g_n(t)$ is symmetric ch.f. we obtain that $g_n$-normalization of $f$ is
\[  {\tilde f}_n(t)= \frac{1}{1-a \ln g_n(t)}.\]
From the condition of c.s.
\[{\tilde f}_n^n(t) = f(t),  \]
and we find that
\[ g_n(t)=\exp\{ \frac{1}{a}(1-(1+a2t^2)^{1/n}\}. \]
\end{ex}

It is known, that Laplace distribution is also geometric stable (so, it is stable for the system with random number of elements) and plays in the case of geometric summation the role of Gaussian distribution. However, we see, it is c.s. for additive systems with non-random number of elements. 

\begin{ex}\label{Ex4}
Let $X$ be a r.v. with gamma distribution. Its Laplace transform is
\[  L(s) = \frac{1}{(1+b s)^{\gamma}},\]
where $b>0$ and $\gamma >0$ are some constants. Introduce normalizing Laplace transform as
\[  g_n(s)=\exp\{ \frac{1}{b}(1-(1+b s)^{1/n}\}.  \]
It is easy to verify, that $X$ is c.s with $g_n$-normalization.
\end{ex}
For $\gamma =1$ we have that exponential distribution is c.s. Usually, geometric distribution is considered as discrete analogue of exponential distribution. It is easy to verify, that geometric distribution with p.g.f.
\[ P(z) =\frac{1-p}{1-p z}, \; p \in (0,1) \]
is c.s. with normalizing p.d.f.
\[ Q_n(z) = \frac{1}{p}(1-(1-p)^{1-1/n}(1-pz)^{1/n}). \]

A part of the next Example represents a hypothesis about p.g.f., and another part gives the statement about c.s. distribution

\begin{ex}\label{Ex5} 
Consider the following p.g.f. 
\begin{equation}\label{eq13}
{\mathcal P}(z) =\exp\{ -\lambda ((1-a z)^\gamma -(1-a)^\gamma \} ,
\end{equation}
where parameters $a$ and $\gamma$ belong to interval $(0,1]$. For $a=1$ we obtain p.g.f. of discrete stable distribution, introduced in \cite{StvHarn79}. Of course, discrete stable distribution is c.s. Our hypothesis consists in statement that the function (\ref{eq13}) is c.s. for all $a, \;\gamma \in (0,1]$. This hypothesis is equivalent to the approval that the normalizing function
\begin{equation}\label{eq14}
Q_n(z) = \frac{1}{a}\Bigl( 1-((1-1/n)(1-a)^{\gamma} +\frac{1}{n}(1-a z)^{\gamma })^{1/\gamma}\Bigr)
\end{equation}
is p.g.f. Unfortunately, we cannot prove this for all values of $\gamma \in (0,1]$. We have proofs for $\gamma=1/2$ and $\gamma=1/3$. Because the proofs for both cases are similar, it supports our hypothesis. Let us give the proof for $\gamma = 1/2$. Really, let us write
\[ Q_n(z) = \sum_{k=0}^{\infty}b_k z^k. \]
The coefficients $b_k$ may be found as
\[ b_0 = \frac{(2-2\sqrt{1-a}+a(n-1))(n-1)}{a n^2}, \; b_1=\frac{1+\sqrt{1-a}(n-1)}{n^2} \]
and
\[  b_k =\frac{2 \sqrt{1-a}a^{k-1}(-1)^{k-1}(n-1) \dbinom{1/2}{k} }{n^2}   \; k=2, 3, \ldots .\]
It is clear, that $b_k \geq 0$ for all non-negative integer $k$. From this we see, that the distribution (\ref{eq13}) is c.s. for $\gamma =1/2$. The same is true for $\gamma =1/3$.
\end{ex}

Let us mention, that c.s. distributions given in Examples \ref{Ex2} - \ref{Ex5} do not have heavy tails. It is opposite with general view on self-similarity of systems and stability. Most popular point of view is that the tails of distributions in self-similar systems are heavy (excluding normal distribution). Now we see, it is not so in the case of random self-similarity (casual stability). However, sometimes the tails of c.s. distributions are heavy. Corresponding examples are, say, classical stable distributions with index $\alpha <2$, limit case of Example \ref{Ex5} as $a \to 1$ with $\gamma <1$. Some examples are given in \cite{KS}, too.

\subsection{Pursuit stability of additive systems}

\begin{ex}\label{Ex6}
Let us consider Sibuya distribution with p.g.f.
\[ {\mathcal P}(z) = 1-(1-z)^\gamma ,\]
where $\gamma \in (0,1)$. For this distribution to be p.c.s. it is necessary and sufficient that
\[  Q_n(z) =1-(1-(1-(1-z)^{\gamma})^n )^{1/\gamma} \]
to be p.g.f. It is possible to prove that $Q_n(z)$ is p.g.f. for any $\gamma \in (0,1)$ and sufficiently large integer $n$ (we do not give the proof here).
\end{ex}

\begin{ex}\label{Ex7}
It is easy to verify that th function
\[ {\mathcal P}(z) = \Bigl (\frac{1-\sqrt{1-z^2}}{z} \Bigr)^M,  \]
is p.g.f. In \cite{KKRT} it was shown that the function
\[ Q_n(z) = 1/T_n(1/z), \]
where $T_n(x)$ is Chebyshev polynomial, is p.g.f. By substitution we see, that
\[ {\mathcal P}(Q_n(z))= {\mathcal P}^n(z). \]
This means that ${\mathcal P}(z)$ is p.g.f. of p.c.s. distribution.
\end{ex}
  
  It is obvious that p.c.s. distribution from Example \ref{Ex7} is not infinitely divisible.
  
\subsection{Multiplicative casual stability}

Now we pass to examples of multiplicative c.s. distributions.

\begin{ex}\label{Ex8}
Let us consider Log-normal distribution with probability density function
\[  p(x) =\exp\{ -\ln^2 x /(2 b^2)\}/(\sqrt{2\pi}b x) \]
and corresponding Mellin transform
\[ M(u) =\exp\{b^2 u^2 /2\} .\]
It is easy to see, that this distribution is multiplicative c.s. with degenerate normalization $1/\sqrt{n}$.
\end{ex}

\begin{ex}\label{Ex9}
Next example is given by double Pareto distribution with the probability density function:
\[ p(x) = \begin{cases} \frac{a2 -1}{2 a}x^a, & x \in (0,1); \\ \frac{a2 -1}{2 a}x^{-a}, & x \geq 1,
\end{cases} \]
where $a>1$ is a parameter. Corresponding Mellin transform has form
\[ M(u) = \frac{a2-1}{a^2-u^2}. \]
Ch.f. of random normalization is
\[ g_n(t)=\exp\{ a-\frac{1}{a}(a^2-1)^{1-1/n}(a^2+t^2)^{1/n}\}. \]
Therefore double Pareto distribution is product c.s.
\end{ex}

Note that the tails of p.c.s. distribution are more heavy than for corresponding additive c.s. distributions. It is easy to see, because for positive random variables multiplicative system may be transformed to additive by passing from r.v.  to its logarithm. In this situation there are some p.c.s. distributions with logarithmic tails.

\begin{ex}\label{Ex10}
Consider a distribution with probability density function
\[  p(x) = \begin{cases} \frac{\exp\{-1/(2 \ln x)\}}{\sqrt{2 \pi}x (\ln x)^{3/2}},& x>1\\ 0, & x \leq 1
\end{cases}\]
Its Mellin transform is
\[ M(u) = \exp\{-\sqrt{2}\sqrt{u}\} \]
for $u>0$. Now it is easy to see, that this distribution is p.c.s. with degenerate normalization. It is analogue of L\'{e}vy distribution, but obviously has logarithmic tail.
\end{ex}

\begin{ex}\label{Ex11}
Let us consider Pareto distribution itself. It has probability density function
\[ p(x)=\frac{\alpha}{x^{1+\alpha}}, \; \alpha >0, \]
for $x \geq 1$. Pareto distribution is geometric product stable (see \cite{KMR}). Here we show, it is also pursuit p.c.s. Really, Mellin transform of Pareto distribution is
\[ M(u) = \frac{\alpha}{\alpha - u}. \]
Ch.f. of random normalization is
\[ g_n(t) = \exp\{ \alpha ((1+it/\alpha )^ n - 1 )\}. \]
\end{ex}

\subsection{Casual stability of min-systems}

Let us pass to examples of  systems with taking minimum operation. 

\begin{ex}\label{Ex12}
Consider c.d.f. of Weibull distribution
\[ F(x) = \exp\{ -\Bigl(\frac{x}{\beta}\Bigr)^{\alpha} \},  \]
for $x>0$, where $\alpha$ and $\beta$ are positive parameters. Define normalization c.d.f. as
\[ G_n(u) = 1-\exp(-a_n u). \]
We have
\[  {\tilde F}_n(x) = {\bar F}(a_n x)= \exp (-a_n^{\alpha}x^{\alpha}/\beta^{\alpha}). \]
From here it is clear, that Weibull distribution is min-c.s. for the case $a_n= 1/n^{1/\alpha}$, and it is pursuit min-c.s. for the case $a_n=n^{1/\alpha}$. Of course, it is natural to call Weibull distribution as min-stable distribution. This terminology was used long ago.
\end{ex}

\begin{ex}\label{Ex13}
Let us show that Gompertz-Makeham distribution is min-p.c.s. Really, its survival function has form
\[  {\bar F}(x) = \exp\{ \xi (1-e^{\lambda x})\}, \;\; \xi >0, \; \;\lambda >0, \]
for positive $x$, and zero for negative values of $x$. From equation (\ref{eq12}) we find the function 
${\bar G}(x)$ for random normalization. It has form
\[ {\bar G}(x) = (1+n(e^{\lambda x}-1)^{-1/\lambda}. \]
Of course, this function is a survival function. Now we see, that Gompertz-Makeham distribution is min-p.c.s.
\end{ex}

\begin{ex}\label{Ex14}
Consider now Pareto distribution. Its survival function is
\[ {\bar F}(x) = 1/x^\alpha; \; \; x \geq 1; \;\; \alpha>0. \]
Normalizing survival function has to be found from equation
\[ 1/(\ln (1/{\bar G}_n(x)))^{\alpha} =1/x^{\alpha /n}. \]
It is 
\[  {\bar G}_n(x) = e^{-x^{1/n}}.\]
Now we see that Pareto distribution is min-c.s. Changing $n$ by $1/n$ we see that Pareto distribution is also min-p.c.s.
\end{ex}

The number of min-c.s. distribution may be prolonged essentially, but we restrict ourself in this paper with  given examples.

\subsection{Casual stability for systems with random number of elements}

We will not consider this model in details. Let us mention, that geometric stable distributions for additive and min systems were studied in many papers and books. The same is true for $\nu$-stable distributions. Therefore, we will give only one non-standard example with max operation.

\begin{ex}\label{Ex15}
Let us consider a family of r.v. $\{\nu_{\alpha}, \; \alpha \in (0,1)\}$ with p.g.f. ${\mathcal P}(z) = 1-(1-z)^{\alpha}$. This is a family of Sibuya distributions. For $\nu_{\alpha}$-max-c.s. distribution we have equation
\[ F(-\ln {\bar G}_a(x)) = 1-(1-F(x))^{\alpha}, \]
where $F$ is a c.d.f., and ${\bar G}_a(x)$ is a survival function. Parameter $a$ depends on $\alpha$.
Let us choose
\[  {\bar G}_a(x) =\exp(-a x). \]
For this situation the function 
\[  F(x) =1- \exp(-\lambda x^b)  \] 
is $nu_{\alpha}$ max-c.s. It is sufficient to put $a=\alpha^{1/n}$.
\end{ex}

\section{A remark on possible connection with fractals}

Analysis of all definitions given above shows that the functions (say, $-\ln g_n(t)$ is Definition \ref{Def1}) used for random normalization are commutative with superposition operation (that is $\ln g_n(\ln g_m(t))=\ln g_m(\ln g_n(t))$). It is interesting to describe all such commutative families of functions. In general statement such problem is very difficult. However, in paper \cite{KKRT} this problem was considered in a very special case of rational p.g.f. for $\nu$-stable distributions. There were used the results obtained by G. Julia and P. Fatou. The methods used by them is connected to dynamical systems on complex plain. G. Julia and P. Fatou had shown, that for commutative rational function corresponding Julia sets are the same. The Julia sets typically have fractal structure. Therefore, there is a hope that similar connection will be present in more general situation. We hope, it will be possible to have more deep connection between fractals and casual stability.

As it was mentioned above, we do not see now any non-trivial connection between random self-similarity and heaviness of tails of corresponding probability distributions.

\section*{Acknowledgments}
Two first authors were supported by the Grant P 203/12/0665  GA\v{C}R.

\end{document}